\documentclass[11pt]{amsart}
\usepackage{amssymb,amsthm,amsmath}
\usepackage[numbers,sort&compress]{natbib}
\usepackage{color}
\usepackage{graphicx}
\usepackage{tikz}
\usepackage[colorlinks,
            linkcolor=blue,
            anchorcolor=blue,
            citecolor=blue
            ]{hyperref}

\newcommand{\SL}{\mathrm{SL}(2,\mathbb{R})}

\def\N{{\mathbb N}}
\def\Z{{\mathbb Z}}
\def\R{{\mathbb R}}
\def\T{{\mathbb T}}

\def\Q{{\mathbb Q}}
\def\tB{\mathcal{B}}
\def\tC{\mathcal{C}}
\def\tF{\mathcal{F}}
\def\tA{\mathcal{A}}

\def\skipaline{\removelastskip\vskip12pt plus 1pt minus 1pt}

\def\Proof{\removelastskip\skipaline
\noindent \it Proof.  \rm}

\newcommand{\la}{\langle }
\newcommand{\ra}{\rangle }

\def\a{\alpha}

\def\ti{\tilde}

\let\newpf\proof \let\proof\relax

\newcommand{\bQ}{\overline{Q}}
\newcommand{\ba}{\overline{A}}

\newcommand{\CD}{\rm CD}

\def\be{\begin{equation}}
\def\ee{\end{equation}}

\def\ba{{\begin{align}}}
\def\ea{{\end{align}}}

\def\bm{\begin{matrix}}
\def\em{\end{matrix}}

\def\a{{\alpha}}

\def\SL{{\mathrm{SL}}}

\def\0{{\mathbf 0}}

\def\cal{\mathcal}

\newtheorem{Theorem}{Theorem}[section]
\newtheorem{Lemma}{Lemma}[section]
\newtheorem{Proposition}{Proposition}[section]
\newtheorem{Corollary}{Corollary}[section]
\newtheorem{Remark}{Remark}[section]

\newtheorem{Definition}{Definition}[section]

\numberwithin{equation}{section}

\def \bn {\hfill \\ \smallskip\noindent}

\theoremstyle{definition}

\def\proof{\bn {\bf Proof.} }

\def\B0{{\bold{0}}}

\catcode`\@=12

\def\Empty{}
\newcommand\oplabel[1]{
  \def\OpArg{#1} \ifx \OpArg\Empty {} \else
    \label{#1}
  \fi}

\newcommand{\comm}[1]{}
\newcommand{\comment}[1]{}

\begin{document}
\setlength{\columnsep}{5pt}

\title[linearization beyond brjuno condition]{Linearization of  quasiperiodically forced circle flow beyond Brjuno condition}

\author{Rapha\"el Krikorian}
\address{
Department of Mathematics, CNRS UMR 8088,
Universit\'{e} de Cergy-Pontoise,  2, av. Adolphe Chauvin F-95302 Cergy-Pontoise, France } 
\email{raphael.krikorian@u-cergy.fr}

\author{Jing Wang}
\address{
Department of Mathematics, Nanjing University of Science and Technology, 210094, China}
\email{jing.wang@njust.edu.cn}

\author {Jiangong You}
\address{
Chern Institute of Mathematics, Nankai University, Tianjin 300071, China
} \email{jyou@nju.edu.cn}

\author{Qi Zhou}
\address{
Department of Mathematics, Nanjing University, Nanjing 210093, China } \email{qizhou@nju.edu.cn}

\date{\today}
\maketitle
\begin{abstract}
We prove that  an analytic quasiperiodically forced circle flow with a
 not super-Liouvillean  base frequency and which is close enough  to some constant rotation is $C^{\infty}$ rotations
reducible, provided its fibered rotation number is Diophantine with respect to the base frequency. As a
corollary, we obtain that among such systems, the linearizable ones and those  displaying
mode-locking are  locally dense for the $C^\infty$-topology.
\end{abstract}

\section{Introduction}

An analytic quasiperiodically forced (qpf) circle flow is the flow   defined
by a  differential equation of the form
\begin{equation}\label{qpf-f}
\left\{
\begin{array}{ll} & \dot{\theta} =
f(\theta,\varphi) \\ & \dot{\varphi} = {\omega}
\end{array} \right.
\end{equation}
where $f: \T^1\times \T^d \rightarrow \R$ is  analytic (here $\T^d=\R^d/\Z^d$), and
${\omega}\in \T^d$ is rationally independent.  We denote by
$({\omega},f)$ the flow of the vector field $(\ref{qpf-f})$. Its discrete version is the qpf
circle diffeomorphism
\begin{equation}\label{qpf-m} (\tilde{\omega},f)=
\tilde{f}: \ (\theta,\varphi)\ \mapsto\ (f_{\varphi}(\theta),
\varphi+ \tilde{\omega})
\end{equation}
where $\theta\mapsto f_{\varphi}(\theta)$ are orientation preserving circle diffeomorphisms and $ (1,\tilde{\omega})\in \T^{d}$ is
 rationally independent.

The qpf circle flow $(\ref{qpf-f})$ and its
discrete version ($\ref{qpf-m}$)
 have been extensively
studied as a source of examples of interesting dynamics and as
models of fundamental phenomena in mathematical physics. One well
known example is the quasiperiodically forced Arnold circle map,
which serves as a simple model for oscillators forced with two or
more incommensurate frequencies \cite{DGO89}. Another example
 comes from the  projective action of  quasi-periodic cocycles:   a typical example is the Harper map, appearing in the study of
quasi-periodic crystals and which can be seen as the projective action of the cocyle associated to an
almost Mathieu operator \cite{AJ05,AK06}.

One of  the classical questions about system $(\omega, f)$  is
whether  there exist a $C^r$ map $H$ : $\T^1 \times \T^d
\rightarrow \T^1\times\T^d$ and $\widetilde{\rho}\in \R$, such that $H$
conjugates the system $({\omega},\rho+f)$ to
$({\omega},\widetilde{\rho})$. If it is the case, we say the system $({\omega},\rho+f)$ is $C^r$ ($r=\infty, \omega$) linearizable, or $C^r$ reducible. 

We first review the corresponding linearization results in the case
$d=1$. In that case,  this means $(\ref{qpf-f})$ is periodically forced  or
$(\ref{qpf-m})$ is unforced, i.e. is a circle map. One of the most
important results for these maps is Poincar\'{e}-Denjoy's
classification theorem: if the rotation number  of an orientation preserving homeomorphism of the circle $f$ is rational, then $f$
has
 periodic orbits; if the rotation number is irrational, then $f$ is {\it
semi-conjugate} to the  corresponding irrational (rigid) rotation;  if the rotation
number is irrational and furthermore $f$ admits a derivative with bounded variation, then the conjugacy is a homeomorphism.
 Subsequently,  Arnold \cite{AR} proved that  if the rotation number is Diophantine,
and if the diffeomorphism  is $C^{\omega}$ close to a rotation, then
it is $C^\omega$ linearizable.   Later in \cite{He79,Yo84}, Arnold's
result was improved obtaining a {\it global} linearization result in the
$C^\infty$ setting: if the rotation number of the diffeomorphism is
Diophantine, then it is linearizable. The $C^\omega$ version was
treated by Yoccoz \cite{Yo02}, in fact, Yoccoz  was able to prove
more: the Brjuno condition is optimal for the local linearization
result, and $\mathcal{H}$ condition \footnote{See \cite{Yo02} for
the definition of $\mathcal{H}$.} is optimal for the global
linearization.

 In this paper, we  focus our attention on the case $d=2$, and without loss of generality, we just assume
$\omega=(1,\alpha)$ with $\alpha  \in
\mathbb{R}\backslash\mathbb{Q}$. From the topological  point of view,
 J\"{a}ger
and Stark \cite{JS06} first gave a Poincar\'{e}-like classification
theorem: if the  quasiperiodically forced circle diffeomorphism
$(\alpha,f)$ is {\it{$\rho$-bounded}}, then either there exists a
{\it{p,q-invariant strip}} and $\rho_{f}$,$\alpha$ and 1 are
rationally dependent or $(\alpha,f)$ is {\it{semi-conjugate}} to the
irrational torus translation
$(\theta,\varphi)\mapsto(\theta+\rho_{f},\varphi+\alpha)$; if
$(\alpha,f)$ is {\it{$\rho$-unbounded}}, then neither of these
alternatives can occur and the map is always topologically
transitive.  Here, $\rho_{f}=\rho(\alpha, f)$ is the {\it fibred rotation
number} of $(\alpha,f)$ (see section \ref{sec:2-3} for details).
However, the analytic part is not  completely satisfactory,  the only result  in that direction being
Herman's local linearization result \cite{He79}, which states that
 the qpf  flow $(\omega, \rho+f)$
\begin{equation}\label{qpf-foced}
\left\{
\begin{array}{ll} & \dot{\theta} =\rho+
f(\theta,\varphi) \\ & \dot{\varphi} = {\omega}
\end{array} \right.
\end{equation}
with  $\rho \in \R$, $f: \T^1\times \T^2 \rightarrow \R$
analytic and $(\omega, \rho_f)$ satisfying the Diophantine
condition
\begin{equation}\label{non-reso-3}
 |\la k,\omega \ra+ l\rho_f|\geq \frac{{\gamma}}{(|k|+|l|)^\tau} \qquad
\textrm{for\ all}\quad (k,l)\in\mathbb{Z}^2\times\mathbb{Z}, \quad |k|+|l| \neq 0,
\end{equation} with ${\gamma}>0,\tau>2$, is  $C^\omega$ linearizable provided  $f$ is sufficiently small (how small depending on $\tau$ and $\gamma$).

 We note that  if \ $l=0$, then $(\ref{non-reso-3})$
implies that $\omega$ is Diophantine  with exponent $\tau$ and constant $\gamma$ which we denote by $\omega\in
DC({\gamma},\tau)$. By Yoccoz's result \cite{Yo02}, Brjuno condition
is optimal for local linearizaton of analytic circle diffeomorphism.  It is
a long-standing question whether
 Brjuno condition  is  also optimal for local  linearization  problems of  other nonlinear systems.
 The essential contribution of this paper  is that   we can really
say something if $\omega$ is not Brjuno or even if it is in some   Liouvillean class. It is clear the
notion of linearization for $(\ref{qpf-foced})$ with Liouvillean
frequencies is too restrictive, since even  $(\omega,g(\varphi))$ is
in general not linearizable if $\omega$ is Liouvillean. The suitable
 notions in this regard are {\it rotations reducibility} and {\it almost
reducibility}.   We say a qpf circle flow $({\omega}, f(\theta,
\varphi))$ is $C^{r} (r= \infty, \omega)$ {\it rotations
reducible} if   it can be $C^{r}$ conjugated to
$({\omega},g(\varphi))$.   We say that $({\omega}, f(\theta, \varphi))$ is
$C^\infty$ {\it almost reducible}, if there exist sequences  $(\rho_n)_{n\in\N}\in \R^\N$,
$(H_n)_{n\in\N}\in (C^\infty(\T\times\T^2,\R\times\R^2))^\N$ and $(f_n)_{n\in\N}\in (C^\infty(\T\times\T^2,\R))^\N$ such that  $(f_{n})$ goes to 0 in the $C^\infty$-topology and $H_n$
conjugates the system $({\omega}, f(\theta, \varphi))$ to
$({\omega},\rho_n+f_n(\theta,\varphi))$. In the analytic category one should be more careful about the choice of topology. 
Let $W_{s,r}(\T \times \T^2)=\{(\theta, \varphi)\in(\T\oplus i\R)\times
(\T^2\oplus i\R^2):|\mathfrak{Im} \theta|< s, |\mathfrak{Im} \varphi|< r\}$. For a
bounded analytic (possibly matrix valued)
 function $f:\T\times\T^2\rightarrow\R$ defined on $ W_{s,r}(\T \times \T^2)$, let
$$ \|f\|_{s,r}=  \sup_{(\theta, \varphi)\in W_{s,r}
}\big|f(\theta,\varphi)\big| .$$  We denote by $C^\omega_{s,r}(\T
\times \T^2,  *)$ the set of all these bounded analytic functions, where $*$ will usually denote $\R$, $\R\times \R^2$. We say that $({\omega}, f(\theta, \varphi))$ is
$C^\omega$ {\it almost reducible}, if there exist $s>0,r>0$ and sequences $\rho_n\in \R$,
$H_n\in C^\omega_{s,r}(\T\times\T^2,\R\times\R^2)$ and $ f_n\in C^\omega_{s,r}(\T\times\T^2,\R)$ such that  $\lim_{n\to\infty}\|f_{n}\|_{s,r}=0$  and $H_n$
conjugates the system $({\omega}, f(\theta, \varphi))$ to
$({\omega},\rho_n+f_n(\theta,\varphi))$. Since this will be important later on, we notice that if one cannot guaranty the convergence of $(f_{n})$ to 0  on {\it fixed} complex strips  of width $s>0$, $r>0$ but instead knows that there exist sequences $(s_{n})_{n\in\N}$, $(r_{n})_{n\in\N}$, $s_{n}>0$, $r_{n}>0$ such that for any $p\in\N$, $\|f_{n}\|_{s_{n},r_{n}}=O((s_{n}r_{n})^p)$, then $(\omega,\rho+f)$ is $C^\infty$ almost reducible (this is a consequence of the classical Cauchy estimates)\footnote{Let us  mention, though we shall not go in that direction in this paper, that another possible notion of almost reducibility would be to consider the convergence of $(f_{n})$ in the analytic inductive limit topology. }.

Apart from the above motivations, our main inspiration comes from
reducibility theory of quasi-periodic  $ SL(2,\R)$ cocycle,
\begin{eqnarray*}\label{cocycle}
(\alpha,A):&\T^{1} \times \R^2 \to \T^{1} \times \R^2\\
\nonumber &(\theta,v) \mapsto (\theta+\alpha,A(\theta) \cdot v).
\end{eqnarray*} The cocycle
$(\alpha,A)$ is called reducible (resp. rotations reducible), if
there exist $B\in C^\omega(2\T, SL(2, \mathbb{R}))$ and $A_* \in
SL(2,\R)$ (resp. $A_* \in C^\omega(\T, SO(2,\R))$) such that $B$ conjugates $(\alpha,A)$ to $(\alpha, A_*)$ which means that $(0,B)\circ (\alpha,A)\circ (0,B)^{-1}=(\alpha,A_{*})$ or equivalently $B(\cdot+\alpha)A(\cdot)B(\cdot)^{-1}=A_{*}(\cdot)$. Similarly, we can define $C^\infty$ (or $C^\omega$) almost reducibility of the cocycles. 
 If $\alpha$ is
Diophantine, the reducibility results are well developed
~\cite{DS75,E92}. In fact, Eliasson \cite{E92} proved that  if $\alpha$ is
Diophantine, then any cocycle which is close to constant is   $C^\infty$ almost reducible.  A breakthrough came recently, based on the   ``cheap
trick" developed in \cite{FK}, where Avila, Fayad and Krikorian~\cite{AFK}
obtained a local positive measures set of $C^\omega$ rotations reducible
result for any base forcing frequency $\alpha \in \mathbb{R}\backslash\mathbb{Q}$. Meanwhile,
for the continuous time  version, by using KAM and  Floquet
theory, Hou and You~\cite{HoY} proved that if the system is  close
to constant ones, then it is  $C^\infty$ almost reducible, and for
full measure set of rotation numbers, the system is $C^\omega$ rotations
reducible.  Readers can consult \cite{YZ1,ZW} for related results.
The  highlight is Avila's global theory of one-frequency $SL(2,\R)$
cocycles \cite{Aglobal},
especially his almost reducibility conjecture (recent solved by him in  \cite{Aac,A2}),
 which says that if the cocycle has subexponentially growth in a fixed analytic strip, then the cocycle is strongly almost reducible.

If a quasi-periodic $ SL(2,\R)$ cocycle is reducible (resp. almost
reducible),  then its corresponding qpf circle diffeomorphism is
linearizable (resp. almost reducible).
 However, the linearization of general
quasiperiodically forced circle flows is more difficult due to their
nonlinear nature. In this paper, we will overcome these difficulties
and generalize the linear result \cite{AFK, HoY} to general nonlinear quasiperiodically forced circle
flows.

To state  our main theorem precisely,  we first introduce some
notations  that will be used in the sequel.  Let $\omega=(1,\alpha)$
with  $\alpha \in  \mathbb{R}\backslash\mathbb{Q}$ and
$\frac{p_n}{q_n}$ be the best convergents of $\alpha$. We denote
\[\widetilde{U}(\alpha):=\sup_{n>0} \frac{\ln\ln q_{n+1}}{\ln q_n}.\]
 If $\widetilde{U}(\alpha)<\infty$, we say  $\alpha$ is {\it not
super-Liouvillean}.
Then our main results can be stated as the following:

\begin{Theorem}\label{mainthm}
Let  $\ti\rho\in\mathbb{R}$, $\gamma>0$, $s>0$, $r>0$, $\tau>2$, $\omega=(1,\alpha)$ with $\alpha\in\R\backslash\Q$ and
$\widetilde{U}(\alpha)<\infty$. Assume that $\rho(\omega,
\tilde{\rho}+ f(\theta,\varphi)) =\rho_f \in DC_\omega(\gamma,\tau)$ in the
sense
\begin{equation}\label{non-reso-1}
|\la k,\omega \ra+ l\rho_f|\geq \frac{\gamma}{(|k|+|l|)^\tau} \quad
\textrm{for \ all\ } (k,l)\in\mathbb{Z}^2\times \mathbb{Z},\quad l \neq 0.
\end{equation}
Then there exists $\varepsilon=\varepsilon(\tau,\gamma,s,r,\widetilde{U})>0$ such that
if $\|f(\theta,\varphi)\|_{s,r} \leq \varepsilon$, the following
holds:
\begin{enumerate}
\item[(a)] The system $(\omega, \tilde{\rho}+f(\theta, \varphi))$  is $C^\infty$
rotations reducible.
\item[(b)] Moreover,  $(\omega, \tilde{\rho}+f(\theta, \varphi))$ is
$C^\infty$ accumulated 
by $C^\omega$ linearizable quasiperiodically forced  circle flow.
\end{enumerate}
\end{Theorem}

\begin{Remark}
As readers will see from the proof, the system is in fact $C^\infty$ almost reducible w.r.t $\theta$, and  $C^\omega$-almost reducible w.r.t $\varphi$. Also note that $C^\infty$
rotations reducibility immediately implies $C^\infty$ almost reducibility.
\end{Remark}

\begin{Remark}
Comparing Theorem \ref{mainthm} with the result of Herman
\cite{He79}, we handle more frequencies $\omega$ including many Liouvillean
ones. Actually, our result applies even to those $\omega$ whose best approximations contain a subsequence of $(q_n)$ such that $q_{n}=O(e^{e^{q_{n-1}}})$.
However, it is still an open question whether the results hold for all
rationally independent $\omega$.
\end{Remark}

These are the results when the fibred rotation number is
non-resonant with the forcing frequency. The reader can consult
\cite{Ji12} for the totally resonant case. The reader can also compare this result with the following recent result of Krikorian \cite{K15},
 who proved that any circle diffeomorphism which  is H\"older conjugated to a rigid rotation is $C^\infty$ almost reducible.

Besides the linearization of a system, another celebrated question
 is whether the rotation number as a function of
the twist parameter is a ``devil's staircase", which means
mode-locking\footnote{Consult section \ref{sec:2-3} for more
information.} is dense in the parameter interval, and for quasiperiodic cocycle this is  equivalent to the fact that
uniform hyperbolicity is dense in the given parameter space.  When
 formulated in the framework of  Schr\"{o}dinger operators, it means the spectrum of the corresponding operator is a
Cantor set.  Yet, existing results are still restricted to almost
Mathieu operators, for which the question is known as the ``Ten
Martini Problem"~\cite{bellissard/simon:1982,puig,AJ05}. Although it
is still an open question whether uniformly hyperbolic is dense in
the category of  $C^\omega$ quasi-periodic $\SL(2,\R)$-cocycles
which are homotopic to the identity,   it is true for  cocycles
which are close to constant ones \cite{Kr}. As a corollary of our
main theorem, the nonlinear version of \cite{Kr} is also true, which
is the following:

\begin{Corollary}\label{dense-mode}
Under the assumptions of Theorem \ref{mainthm},  the qpf flow $(\omega,
\tilde{\rho}+f(\theta,\varphi))$ is $C^\infty$ accumulated 
 by
mode-locked quasiperiodically forced circle flows.
\end{Corollary}

\begin{Remark}
To the best knowledge of the authors, this result  is the first result of denseness about mode-locking in the nonlinear setting.
   \end{Remark}

\begin{Remark}
All these results are also valid for the discrete case, i.e.,
quasiperiodically forced circle maps.
\end{Remark}

Finally, we want to point out that in this paper  the method of  diagonally
dominant operators we developed  is  new in itself (consult more discussions in section \ref{sec:3}).
The method has recently been proved to be very powerful: it can  be used to solve several interesting problems, and we will
come back to theses issues in forthcoming papers.

\section{Notations and preliminaries}

\subsection{Continued fraction expansion.}\label{sec:2.1}
Let $\alpha \in (0,1)$ be irrational. Define $ a_0=0,
\alpha_{0}=\alpha,$ and inductively for $k\geq 1$,
$$a_k=[\alpha_{k-1}^{-1}],\qquad \alpha_k=G(\alpha_{k-1})=\alpha_{k-1}^{-1}-a_k=\{ \alpha_{k-1}^{-1}\}.$$
We define
$$p_0=0, \qquad p_1=1,$$ $$
q_0=1, \qquad q_1=a_1,$$ and inductively,
\begin{eqnarray*}
 &&p_k=a_kp_{k-1}+p_{k-2},\\&& q_k=a_kq_{k-1}+q_{k-2}.
\end{eqnarray*}
Then  the sequence $(q_n)$  is  the sequence of denominators of the
best rational approximations for $\alpha \in \R \backslash \Q$, since
it satisfies
\begin{equation}
\forall 1 \leq k < q_n,\quad \|k\alpha\|_{\T} \geq
\|q_{n-1}\alpha\|_{\T},
\end{equation}
and
\begin{equation}
\frac{1}{q_n+q_{n+1}}\leq\|q_n \alpha \|_{\T} \leq {1 \over
q_{n+1}},
\end{equation}
where we use the notation
$$\|x\|_{\T}=\inf_{p\in\Z}|x-p|.$$

\subsection{CD bridge.}
For  any  $\alpha \in \R \backslash \Q$, we will fix in the sequel a
particular subsequence ${(q_{n_k})}$ of the denominators of the
continued fraction expansion for $\alpha$,  which is denoted  by
$(Q_k)$, and the subsequence $(q_{n_k+1})$ is denoted by $(\bQ_k)$.

\begin{Definition}\cite{AFK}\label{CDbridge}
Let $0<\tA \leq \tB\leq \tC$. We say that the pair of denominators $(q_l,q_n)$ ($l<n$) forms a
$\CD(\tA,\tB,\tC)$ bridge if
\begin{itemize}
\item $q_{i+1}\leq q_i^{ \tA}, \quad \forall i=l,\ldots,n-1$
\item $q_l^{\tC}\geq q_n\geq q_l^{\tB}$.
\end{itemize}
\end{Definition}

\begin{Lemma}\cite{AFK}\label{CDbridge}
For any ${\tA}>0$, there exists  a subsequence $(Q_k)$ such that $Q_0=1$ and for each $k\geq 0$,
$Q_{k+1}\leq \bQ_k^{{\cal A}^4}$. Furthermore, either $\bQ_k\geq
Q_k^{\cal A}$, or the pairs $(\bQ_{k-1},Q_{k})$ and $(Q_k,Q_{k+1})$
are both $\CD({\cal A},{\cal A},{\cal A}^3)$ bridges.
\end{Lemma}

In the sequel, we assume $\mathcal{A}=8$, and $(Q_n)$ is the
selected subsequence in Lemma \ref{CDbridge} accordingly. As an
immediate corollary of Lemma \ref{CDbridge}, we have

\begin{Corollary}\label{cor2-1}
If $\widetilde{U}(\alpha)<\infty,$ then we have $Q_n\geq
Q_{n-1}^\cal A$ for every $n\geq 1$. Furthermore, we have
$$\sup_{n>0} \frac{\ln\ln Q_{n+1}}{\ln Q_n}\leq U(\alpha),$$ where
$U(\alpha):=\widetilde{U}(\alpha)+4\frac{\ln \mathcal A}{\ln 2}<\infty$.
\end{Corollary}

\Proof For $n=1$, it is obvious that $Q_1\geq Q_0^\mathcal{A}$ since
$Q_0=1$. For $n\geq2$, we distinguish two cases below. If
$\bQ_{n-1}\geq Q_{n-1}^{\cal A}$, then $Q_n\geq \bQ_{n-1}\geq
Q_{n-1}^\cal A$. Otherwise, the pairs $(\bQ_{n-2},Q_{n-1})$ and
$(Q_{n-1},Q_n)$ are both $\CD({\cal A},{\cal A},{\cal A}^3)$
bridges. Thus, we get $Q_n\geq Q_{n-1}^\cal A$ by definition.

For the second statement, we have \[ \frac{\ln\ln Q_{n+1}}{\ln
Q_n}\leq \frac{\ln {\cal A}^4 + \ln\ln \overline{Q}_{n} }{\ln
Q_n}\leq \widetilde{U}+\frac{4\ln\mathcal{A}}{\ln2}=U ,\] since for
each $n\geq0$, $Q_{n+1}\leq \bQ_n^{{\cal A}^4}$.\qed

\subsection{The fibred rotation number and mode-locking.}\label{sec:2-3}
Suppose $({\omega},f)$ is a qpf circle flow. Let
$$\rho(\omega,f)=\lim_{t\to\infty} \frac{\hat{\Phi}_{\varphi}^t(\hat{\theta})}{t}$$
be the fibred rotation number associated with $({\omega},f)$, where
$\hat \Phi_\varphi^t(\hat \theta): \R_+^1\times \R^1\times
\T^2\rightarrow \R^1$, via $(t,\hat \theta, \varphi)\mapsto \hat
\Phi_\varphi^t(\hat \theta) $ denotes the lift of the flow  of
$({\omega},f)$ of the first variable. The limit exists and is
independent of $(\hat \theta,\varphi)$ \cite{He83}. As a direct
consequence of the definition, we have the following well-known results:

\begin{Lemma}\cite{He83} \label{fib-rota-1}
Let $\omega=(1,\alpha)$,   $\tilde{\rho}\in\R$. If
$\|f(\theta,\varphi)\|_{C^0} \leq\varepsilon$, then
$$|\rho(\omega, \tilde{\rho}+f(\theta, \varphi)) - \rho(\omega, \tilde{\rho}) |\leq \varepsilon$$
\end{Lemma}

\begin{Lemma}\cite{He83}\label{fib-rota-2}
Suppose $(\omega, f)$ is a qpf circle flow and $H\in C^0(\T\times
\T^2, \T\times \T^2)$ is homotopic to the identity. Then the fibred
rotation number of $(\omega, f)$ remains the same under the
transformation $H$.
\end{Lemma}

 If the
rotation number of a qpf circle flow $(\omega,f)$  remains constant under all sufficiently small ${\cal
  C}^0$-perturbations of $f$,  we say $(\omega,f)$ is  {\it mode-locked}. For our
  purpose, we need the following proposition, which was proved in
  \cite{BJ}.

\begin{Proposition}\cite{BJ} \label{mode}
Suppose the quasiperiodically forced circle flow $(\omega,f)$ is
given by the projective action of a quasi-periodic $sl(2,\R)$-flow.
Then $(\omega,f)$ is mode-locked if and only if the quasi-periodic
$sl(2,\R)$-flow is uniformly hyperbolic.
\end{Proposition}

Readers can consult \cite{BJ} for more interesting properties of
mode-locking.

\section{Outline of the proof}\label{sec:3}
In this section, we will outline the main ideas of the proof,
and show the principal difference between our method and previous
works. We first introduce some notations that will be used in the
sequel. If $f \in C^\omega(\T^2,\R)$, we denote $
\widehat{f}(k)=\int_{\T^2}f(\varphi)e^{- i \la
k,\varphi\ra}d\varphi$.  For  $f\in C^\omega(\T \times \T^2,\R)$, we
define its fourier coefficients by $f_l^k=\int_{\T \times
\T^2}f(\theta,\varphi)e^{-i l \theta}e^{- i \la k,\varphi\ra}d\theta
d\varphi$. For any $N>0$, we denote the truncation and projection
operators $\mathcal{T}_N$ and $\mathcal{R}_N$ by
$$ \mathcal{T}_N (f)=\sum_{0<|k|+|l|<N} f_l^k  e^{i l\theta}e^{i \la k,\varphi\ra},
\quad \mathcal{R}_N (f) =\sum_{ |k|+|l|\geq N} f_l^k  e^{i
l\theta}e^{i \la k,\varphi\ra}.$$

 The proof  is based on a modified KAM scheme. Considering  the quasiperiodically forced circle flow
 $(\omega,\rho+f(\theta,\varphi))$, if we want to eliminate the non-resonant terms of $f(\theta,\varphi)$ as in usual KAM
 steps, then the classical homological equation on the fibre reads as
$$f(\theta,\varphi)-\partial_{\omega}h-\rho\frac{\partial h}{\partial
\theta}=0.$$ Checking the fourier coefficients, we have
$$
h_{l}^{k}=\frac{f_{l}^{k}}{i\big(\la k,\omega \ra +l\rho \big)}.
$$
Since we have no Diophantine condition on $\omega$, the homological
equation may have no analytic solution. This is the essential
difference compared to the classical KAM theorem, which also means
the resonant terms $f_{0}^{k}e^{i\langle k,\varphi\rangle}$ cannot
be solved at all. So we rewrite the system as $(\omega, \rho
+g(\varphi) + f(\theta,\varphi))$, with $\int_{\T^2}f(\theta,\varphi)d\varphi=0$. In this case, the homological
equation on the fibre is
\begin{equation}\label{coho}
\partial_{\omega}h+(\rho
+g(\varphi))\frac{\partial h}{\partial \theta}=f(\theta,\varphi).
\end{equation}

In order to get desired result, we distinguish three steps. The first step is to eliminate the
non-resonant terms of $g(\varphi)$ by  solving $\partial_\omega
h(\varphi)=\mathcal {T}_{Q_n} g(\varphi)$. Notice that
$\|h(\varphi)\|$ may be very large. However, in this step, we use a
small trick to control $\|\mathfrak{Im}h(\varphi)\|$ at the cost of
reducing the analytic radius greatly.

The second step is to solve the homological equation
$$
\partial_{\omega}h+(\rho+
\widetilde{g}(\varphi))\frac{\partial h}{\partial
\theta}=f(\theta,\varphi),
$$
where $\|\widetilde{g}(\varphi)\|=O(\|f\|).$ By  introducing diagonally dominant operators,  we can solve the approximate equation
$$
\partial_{\omega}h+\rho\frac{\partial h}{\partial \theta}+
\mathcal {T}_K\big(\widetilde{g}(\varphi)\frac{\partial h}{\partial
\theta}\big)=\mathcal {T}_Kf(\theta,\varphi)
$$
 and then  make the perturbation $f(\theta,\varphi)$ as small as we can
by iteration.

Using these two steps, we can already prove our almost  reducibility result.
However,  to obtain the rotations reducibility result,  in the end of one
KAM step  we need to inverse the first step;  the conjugation we get is then
close to the identity.

\section{The inductive step}
\subsection{Basic proposition.}
In this subsection, we will show how to use the method of   diagonally
dominant operators to solve the homological equation
\begin{equation}\label{homo}
\partial_{\omega}h+(\rho
+g(\varphi))\frac{\partial h}{\partial \theta}=f(\theta,\varphi).
\end{equation}
We should point out that the following proposition holds, irrespective of
any arithmetical property on the base frequency $\omega$, and this
is fundamental in our reduction.

\begin{Proposition}\label{Diagonal}
Let $\gamma>0$, $\tau>2$, $r>\sigma>0$, $s>\delta>0$,
$\sigma\leq\delta/2$,    $\rho\in DC_\omega(\gamma,\tau)$,
$g(\varphi)\in C^\omega_r(\T^2,\R)$, $f(\theta,\varphi)\in
C_{s,r}^\omega(\T\times\T^2,\R)$ with
$\int_{\T}f(\theta,\varphi)d\theta=0$. There exist
$0<\widetilde\eta\leq \eta\ll 1$,  such that if
$\|g(\varphi)\|_{r}\leq \eta$, $\|f(\theta,\varphi)\|_{s,r}\leq
\widetilde{\eta}$ and
\begin{equation}\label{dia-2}
K=\left[\frac{1}{\sigma}\ln\frac{1}{\widetilde{\eta}}\right]<
 (\frac{\gamma^{2}}{\eta})^{\frac{1}{2\tau+3}},
\end{equation}
then the homological equation $(\ref{homo})$
 has an approximate solution
$h(\theta,\varphi)$ with estimate
$$\|h\|_{s-\delta,r-\sigma} \leq
\frac{2\widetilde{\eta}}{\gamma\sigma^{2+\tau}},$$ and the error
term $\widetilde{P}=\mathcal{R}_{K}(-g(\varphi)\frac{\partial
h}{\partial \theta}+f(\theta,\varphi))$ satisfies
$$\|\tilde{P}\|_{s-\delta,
r-\sigma}\leq\frac{4\widetilde{\eta}^{2}}{\gamma\sigma^{3+\tau}}.$$
\end{Proposition}

\Proof Since the homological equation $(\ref{homo})$ may  have no
analytic solution, we solve its approximate equation
\begin{equation}\label{coho-equ-5-2}
\partial_{\omega}h+\rho\frac{\partial h}{\partial \theta}
+\mathcal{T}_{K}\big(g(\varphi)\frac{\partial h}{\partial
\theta}\big)=\mathcal{T}_{K}f(\theta,\varphi).
\end{equation}
Let
$$ f(\theta,\varphi)=\sum_{l}f_{l}(\varphi)e^{il\theta},\qquad
f_{l}(\varphi)=\sum_{k}f_{l}^{k}e^{i\langle k,\varphi\rangle},$$
\begin{displaymath}
h(\theta,\varphi)=\sum_{0<|l|<K}h_{l}(\varphi)e^{il\theta},\qquad
h_{l}(\varphi)=\sum_{|k|<K-|l|}h_{l}^{k}e^{i\langle
k,\varphi\rangle}.
\end{displaymath}
In order to solve $(\ref{coho-equ-5-2}),$ it is equivalent to solve
\begin{equation}\label{mat}
\partial_{\omega}h_{l}(\varphi)+il\rho h_{l}(\varphi)
+il\Gamma_{K-|l|}\left(g(\varphi)h_{l}(\varphi)\right)=\Gamma_{K-|l|}f_{l}(\varphi)
\end{equation}
for $0<|l|<K$, where for $ f(\varphi)\in C^\omega(\T^2,\R), N>0,$
$$\Gamma_N f(\varphi)=\sum_{0\leq |k|<N}\widehat{f}(k)e^{i\langle k,\varphi\rangle}.$$

 For any fixed $l$, $(\ref{mat})$ can be viewed  as a matrix equation
$$(A_{l}+G_{l})\overline{h}_{l} = \overline{f}_{l},$$
where
$$\overline{h}_{l}=(h_{l}^{k})^{T}_{|k|< K-|l|}, \qquad  A_l=diag(i\langle k,\omega\rangle +il \rho:|k|< K-|l|),$$
$$\overline{f}_{l}=(f_{l}^{k})^{T}_{|k|< K-|l|}, \qquad  G_{l} = ( i l \widehat{g}(p-q))_{|p|,|q|< K-|l|}.$$
If we denote $\Omega_{l,r'}=diag(\cdots, e^{|k|r'}, \cdots)_{|k|<
K-|l|}$ for any $r'\leq r$, then
$$\Omega_{l,r'}({A}_{l}+{G}_{l})\Omega_{l,r'}^{-1}\Omega_{l,r'}\overline{h}_{l} = \Omega_{l,r'}\overline{f}_l.$$
Rewrite it as
$$(\widetilde{A}_{l,r'}+\widetilde{G}_{l,r'})\widetilde{h}_{l,r'} =
\widetilde{f}_{l,r'},$$ where
$$\widetilde{A}_{l,r'}=\Omega_{l,r'}{A}_{l}\Omega_{l,r'}^{-1}, \qquad \widetilde{G}_{l,r'}=\Omega_{l,r'}{G}_{l}\Omega_{l,r'}^{-1},$$
$$\widetilde{h}_{l,r'}=\Omega_{l,r'}\overline{h}_{l}, \qquad \widetilde{f}_{l,r'}=\Omega_{l,r'}\overline{f}_l, $$ with estimate
$$\|\widetilde{G}_{l,r'}\|\leq|l|(K-|l|)^2\|g\|_r.$$

Since $\rho\in DC_\omega(\gamma,\tau)$, then
together with (\ref{dia-2}), we have
\begin{equation}\label{dio-c-3}
|\langle k,\omega\rangle+l\rho|>(|l|(K-|l|)^2\|g\|_{r})^{1/2}
\end{equation}
for all $0<|k|+|l|< K, l\neq 0$. As a result,  the diagonally
dominant  operators $\widetilde{A}_{l,r'}+\widetilde{G}_{l,r'}$ has
a bounded inverse and $\|(I+\widetilde A_{l,r'}^{-1}\widetilde
G_{l,r'})^{-1}\|_{op( l^{1})}<2$,  where $\|\cdot\|_{op(l^1)}$ denotes the operator norm associated to the $l^1$ norm $|u|_{l^1}=\sum_{|k|<K-| l |} | u^k|$ (indeed, if $M$ is a matrix, $\|M\|_{op(l^1)}\leq C(M)$ with $C(M)=\sum_{i}\max_{j}|M_{ij}|$) and from the above estimates it is clear that $C(\widetilde A_{l,r'}^{-1}\widetilde G_{l,r'})<1/2$.

If   $\sigma\leq\delta/2$, then we can  estimate
\begin{eqnarray*}
&&\|h(\theta,\varphi)\|_{s-\delta/2,
r-\sigma}\\ &\leq&\sum_{|k|+|l|<K}|h_l^k|e^{|k|(r-\sigma)}e^{|l|(s-\delta/2)}=\sum_{0<|l|<K}|\widetilde
h_{l,r-\sigma}|_{l^1}e^{|l|(s-\delta/2)}\\
&=&\sum_{0<|l|<K}|(I+\widetilde A_{l,r-\sigma}^{-1}\widetilde
G_{l,r-\sigma})^{-1}\widetilde A_{l,r-\sigma}^{-1}\widetilde
f_{l,r-\sigma}|_{l^1}e^{|l|(s-\delta/2)}\\
&\leq&\sum_{0<|l|<K}\|(I+\widetilde A_{l,r-\sigma}^{-1}\widetilde
G_{l,r-\sigma})^{-1}\|_{op(l^1)}\cdot|\widetilde A_{l,r-\sigma}^{-1}\widetilde{f}_{l,r-\sigma}|_{l^1}e^{|l|(s-\delta/2)}\\
&\leq&\sum_{0<|l|<K}\sum_{|k|<K-|l|}\frac{2(|k|+|l|)^{\tau}}{\gamma}|f_l^k|e^{|k|(r-\sigma)}e^{|l|(s-\delta/2)}\\
&\leq&\frac{2\widetilde \eta}{\gamma \sigma^{2+\tau}}.
\end{eqnarray*}
Consequently by Cauchy estimates, we get the control of the error
term:
\begin{eqnarray*}
\|\tilde{P}\|_{s-\delta, r-\sigma}&\leq
&\|\mathcal{R}_{K}f(\theta,\varphi)\|_{s-\delta,r-\sigma}+\left\|\mathcal{R}_{K}\left(g(\varphi)\frac{\partial
h}{\partial \theta}\right)\right\|_{s-\delta,r-\sigma}\\
&\leq&
\frac{\widetilde{\eta}e^{-K\sigma}}{\sigma^{2}}+\frac{\eta\|h(\theta,\varphi)\|_{s-\delta/2,r-\sigma}e^{-K\sigma}}{\sigma}\\
&\leq&
\frac{\widetilde{\eta}^2}{\sigma^2}+\frac{2\widetilde{\eta}^{2}}{\gamma\sigma^{3+\tau}}
<\frac{4\widetilde{\eta}^{2}}{\gamma\sigma^{3+\tau}}.
\end{eqnarray*}

 \qed

\subsection{KAM step.}\label{sec-3}

We give details about one step of the KAM iteration. For simplicity,
we introduce some notations which will be used in this section.  For any
$r,s,\eta,\tilde{\eta}>0$, $\rho_f\in\R$, we define
\begin{eqnarray*}
&& \tF_{s,r}(\rho_f, \eta, \tilde{\eta})\\
&=& \left\{\tilde{\rho}+ g(\varphi)+f(\theta,\varphi)\in C^\omega_{s,r}(\T
\times \T^2,\R)  \Bigg{|}   \begin{array}{cc} \rho(\omega, \tilde{\rho}+
g(\varphi)+f(\theta,\varphi))=\rho_f \\  \|g(\varphi)\|_r \leq \eta,
\|f(\theta,\varphi)\|_{s,r}\leq   \tilde{\eta} \end{array} \right\}.
\end{eqnarray*}

Let $\alpha\in \mathbb{R}\backslash\mathbb{Q}$, with
$U=\tilde U(\alpha)+12<\infty$, $(Q_n)$ is the selected sequence of $\alpha$
by  Lemma \ref{CDbridge} with $\mathcal{A}=8$. Let $r_0, s_0, \gamma>0, \tau>2$ and $Q_*$ be the smallest $Q\in\N_0$ such that
\[\ln Q<\frac{Q^{1/4}r_0}{40c\tau U}.\]
Suppose that
$\varepsilon_0$ is small enough such that
\begin{equation}\label{enseq}
\varepsilon_0<\min\{\frac{(r_0s_0\gamma)^{12(\tau+3)}}{\tau!Q_{1}^{2c\tau
U}},e^{-2c\tau U}, e^{-40(\ln Q_*)^2c\tau U}\}\quad\textrm{and}\quad
\ln\frac{1}{\varepsilon_0}<(\frac{1}{\varepsilon_0})^{\frac{1}{12(2\tau+3)}}.
\end{equation}
where $c$  is a global constant with $c>10(\tau+3)/\tau$. For any
given $r_0,s_0,\varepsilon_0$, we inductively define some sequences
depending on $r_0,s_0,\varepsilon_0$ for $j\geq 1$:
$$\Delta_{1}=\frac{s_{0}}{10},\quad
\Delta_{j}=\frac{\Delta_{1}}{2^{j-1}},$$
\begin{equation}\label{seque}
r_{j}=\frac{r_0}{4Q_{j}^3},\quad s_{j}=s_{j-1}-\Delta_{j},
\end{equation}
$$\varepsilon_{j}=\frac{\varepsilon_{j-1}}{Q_{j+1}^{2^{j+1}c\tau
U}},\quad \widetilde\varepsilon_j=\sum_{m=0}^{j-1}\varepsilon_m,$$
$$K^{(j)}=\left[\Big(\frac{\gamma^{2}}{4\varepsilon_{j}}\Big)^{\frac{1}{2\tau+3}}\right].$$
Since $\frac{1}{Q_j^3}$ goes to 0 much faster than $\Delta_j$, we
can just assume that $4r_j<\Delta_j$ without loss of generality.

\subsubsection{Eliminate the non-resonant terms.}\label{sec-4-2-1}

\begin{Lemma}\label{lemma4-1}
Given qpf circle flows
\begin{eqnarray} \label{sys-13}\left\{
\begin{array}{ll} & \dot{\theta} =
\rho_f+g(\varphi)+f(\theta,\varphi)
\\ & \dot{\varphi} = \omega=(1,\alpha)
\end{array} \right.
\end{eqnarray}
if  $\rho_f+g(\varphi)+f(\theta,\varphi)\in
\tF_{s_{n-1},r_{n-1}}(\rho_f,
4\widetilde\varepsilon_{n-1},\varepsilon_{n-1})$, then, if we denote
  $\overline{s}_{n}=s_{n-1}- \frac{\Delta_n}{3},$
$\overline{r}_n=\frac{r_0}{Q_n^3}$,  there exists
$h(\varphi)$ with
$\|h(\varphi)\|_{\overline{r}_n}\leq
Q_n^{\frac{7}{4}}\varepsilon_0^{\frac{1}{2}},$ such that the
transformation $\theta = \overline{\theta}+h(\varphi)$ (mod 1)
conjugates the system $(\ref{sys-13})$  into
\begin{eqnarray}\label{sys-2} \left\{
\begin{array}{ll} & \dot{\overline{\theta}} =
\rho_{f}+\overline{g}(\varphi)+\overline{f}(\overline{\theta},\varphi)
\\ & \dot{\varphi} = \omega=(1,\alpha)
\end{array} \right.
\end{eqnarray}
with $
\rho_{f}+\overline{g}(\varphi)+\overline{f}(\overline{\theta},\varphi)\in
\tF_{\overline{s}_{n},\overline{r}_n  } (\rho_f,
\varepsilon_{n-1}^{1/2},\varepsilon_{n-1}).$
\end{Lemma}
\Proof Under the transformation $ \theta =
\overline{\theta}+h(\varphi)$ (mod 1), the fibred equation becomes
$$\dot{\overline{\theta}}=\rho_f-\partial_{\omega}h+g(\varphi)
+f(\overline{\theta}+h(\varphi),\varphi).$$ Let
$\partial_{\omega}h(\varphi)=\mathcal{T}_{Q_n}g(\varphi)$, then the fibred
equation is
$$\dot{\overline\theta}=\tilde{f}(\overline{\theta},\varphi)=\rho_f+\widehat g(0)+\mathcal R_{Q_n} g(\varphi)+f(\overline{\theta}+h(\varphi),\varphi).$$
Since for $0<|k|<Q_n$, $|\langle k,\omega\rangle|>\frac{1}{2Q_n}$,
it is clear that
\begin{equation}
\|h(\varphi)\|_{\frac{r_{n-1}}{2}}\leq
2Q_n\sum_{0<|k|<Q_n}\|g\|_{r_{n-1}}e^{-|k|\frac{r_{n-1}}{2}}
\leq\frac{64Q_n\varepsilon_0}{r_{n-1}^2}\leq
Q_n^{\frac{7}{4}}\varepsilon_0^{\frac{1}{2}}.
\end{equation}
Observe that if  $\varepsilon_0^{-1/2}<Q_n^{7/4}$,  we may
lose the control of the norm of
$f(\overline{\theta}+h(\varphi),\varphi)$.

However,   in order to estimate
$f(\overline{\theta}+h(\varphi),\varphi)$, it is sufficient to
control the imaginary part of $h(\varphi)$.  To fulfill this, we
need a small trick,  which says that $\|\mathfrak{Im}h(\varphi)\|$
can be  well controlled  at the cost of reducing the analytic radius
greatly.

Let $\varphi=\varphi_{1}+ i\varphi_{2}$ for $
\varphi_{1}\in \mathbb{T}^{2}$, $\varphi_{2}\in \mathbb{R}^{2}$ and define
$$h_{1}(\varphi_{1})=\sum_{0<|k|<Q_n}\frac{\widehat{g}(k)}{i\langle k,\omega\rangle}e^{i\langle k,\varphi_{1}\rangle}, \quad\quad
h_{2}(\varphi)=h(\varphi)-h_{1}(\varphi_{1}).$$ Since $g(\varphi)$
is real analytic, we have $\mathfrak{Im}h_1(\varphi_1)=0$. Moreover,
by Corollary \ref{cor2-1}, we have $Q_n\geq Q_{n-1}^8$, which
implies the following estimate
\begin{eqnarray*}
\|\mathfrak{Im}h(\varphi)\|_{\overline{r}_n}&=&\|\mathfrak{Im}h_{2}(\varphi)\|_{\overline{r}_n}\leq \|h_{2}(\varphi)\|_{\overline{r}_n}\\
&\leq&2Q_n\sum_{0<|k|<Q_n}|\widehat{g}(k)||e^{-\langle k,\varphi_{2}\rangle}-1|\\
&<&2Q_n\sum_{0<|k|<Q_n}\|g\|_{r_{n-1}} e^{-|k|(r_{n-1}-\bar r_n)}\cdot|k|{\overline{r}_n} \\
&<& \frac{C\varepsilon_0r_0}{Q_n^{2}(r_{n-1}-\overline{r}_n)^{3}} <
\frac{\varepsilon_0^{1/2}}{Q_n^{1/2}} <\frac{\Delta_{n}}{3}.
\end{eqnarray*}
As a consequence, we have
$$\|f(\overline{\theta}+h(\varphi),\varphi)\|
_{\overline{s}_n,\overline{r}_n}\leq
\|f(\theta,\varphi)\|_{s_{n-1},r_{n-1}}\leq \varepsilon_{n-1},$$
$$
\| \cal R_{Q_n} g(\varphi)\|_{\overline{r}_n}\leq CQ_n^2\varepsilon_0
e^{-Q_nr_{n-1}/2} \leq\varepsilon_{n-1}^{1/2}/3, $$ by (\ref{enseq}).
Lemma \ref{fib-rota-2} implies  $\rho(\omega, \tilde{f}(
\overline{\theta},\varphi ))=\rho_f$. Then by Lemma
\ref{fib-rota-1}, it follows that  $$|\widehat g(0)|\leq
\|f(\overline{\theta}+h(\varphi),\varphi)\|
_{\overline{s}_n,\overline{r}_n}+\| \cal R_{Q_n}
g(\varphi)\|_{\overline{r}_n}\leq 2\varepsilon_{n-1}^{1/2}/3.$$ Let
$\overline g(\varphi)=\widehat g(0)+\cal R_{Q_n}g(\varphi)$ and
$\overline f(\overline \theta,\varphi)= f(\overline
\theta+h(\varphi), \varphi)$. Then the result follows. \qed

\subsubsection{Reduction by diagonally dominant operators}\label{sec-4-2-2}
In this section, we will apply Proposition \ref{Diagonal} to make
the perturbation as small as we can.

\begin{Lemma}\label{lemma4-2}
 Under the
assumptions of Lemma \ref{lemma4-1}, if furthermore $\rho_f\in
DC_\omega(\gamma,\tau)$,
 then there exists $\overline{H}\in
C_{s_{n+},r_{n+}}^\omega(\T\times\T^2,\T\times\T^2)$ with estimates
\begin{equation}\label{tran-1}
\|\overline{H}-id\|_{s_{n+},r_{n+}}\leq
4\varepsilon_{n-1}^{\frac{3}{4}},
\end{equation}
\begin{equation}\label{tran}\|D(\overline{H}-id)\|_{s_{n+},r_{n+}}\leq
4\varepsilon_{n-1}^{\frac{3}{4}},\end{equation} such that
$\overline{H}$ conjugates the system $(\ref{sys-2})$ to
\begin{eqnarray}\label{sys-12} \left\{
\begin{array}{ll} & \dot{\bar\theta}_+ = \rho_{f} +\overline{g}_{+}(\varphi) +
\overline{f}_{+}(\bar\theta_+,\varphi) \\ & \dot{\varphi} =
\omega=(1,\alpha)
\end{array} \right.
\end{eqnarray}
with $ \rho_{f} +\overline{g}_{+}(\varphi) +
\overline{f}_{+}(\bar\theta_+,\varphi)\in \tF_{s_{n+},r_{n+}}(\rho_f,
2\varepsilon_{n-1}^{1/2},\varepsilon_{n} )$ and $\|\overline g_+-\overline g\|_{r_{n+}}\leq 4\varepsilon_{n-1}$, where we denote
$r_{n+}=\frac{r_0}{2Q_n^3}$,
$s_{n+}=\overline{s}_{n}-\frac{\Delta_{n}}{3}$.

\end{Lemma}
\Proof To avoid ambiguous notations, in the proof of this lemma, we
fix $n$, and denote temporarily
$\widetilde{r}=\overline{r}_n=\frac{r_0}{Q_n^3},$
$\widetilde{s}=\overline{s}_n=s_{n-1}-\frac{\Delta_n}{3}$,
$\eta=2\varepsilon_{n-1}$, $\tilde \eta=2\varepsilon_{n-1}^{1/2}$. We will prove this lemma by iteration.
First we define sequences:
$$ \widetilde{r}_{0}=\widetilde{r},\qquad \widetilde{s}_{0}=\widetilde{s},\qquad \eta_{\nu}=\eta^{(\frac{3}{2})^{\nu}},$$
$$\sigma_{1}=\frac{\widetilde{r}}{4}=\frac{r_0}{4Q_n^3},\qquad
\sigma_{\nu+1}=\frac{1}{2^{\nu}}\sigma_{1},$$$$\delta_{1}=\frac{\Delta_{n}}{6},\qquad
\delta_{\nu+1}=\frac{1}{2^{\nu}}\delta_{1},$$
$$\widetilde{r}_{\nu}=\widetilde{r}_{\nu-1}-\sigma_{\nu},\quad \widetilde{s}_{\nu}=\widetilde{s}_{\nu-1}-\delta_{\nu}.$$
$$K_{\nu}=\left[\frac{1}{\sigma_{\nu}}\ln\frac{1}{\eta_{\nu-1}}\right].$$
Since $\frac{1}{Q_n^3}$ goes to 0 much faster than $\Delta_n$, we
can assume that $\widetilde{r}<\Delta_n/3$ without loss of
generality. Then as a consequence $\sigma_{\nu}\leq\delta_{\nu}/2$.

Let $N=\left[2^nc_1\tau U\ln Q_n\right]+1$, where
$c_1=\frac{c}{32(\tau+3)\ln 3}$. Assume for $j=1,2,\cdots, \nu-1\ <\
N$, there are $h_j\in C_{\tilde s_{j-1}-\frac{\delta_j}{2}, \tilde
r_{j-1}-\frac{\sigma_j}{2}}^\omega(\T\times\T^2, \R), \bar f_j\in
C_{\tilde s_j, \tilde r_j}^\omega(\T\times\T^2,\R), \bar g_j\in
C_{\tilde r_{j-1}}^\omega(\T^2, \R)$, with 
$\|\bar
f\|_{\tilde s_j,\tilde r_j}\leq \eta_j$, $\|h_j\|_{\tilde
s_{j-1}-\frac{\delta_j}{2}, \tilde r_{j-1}-\frac{\sigma_j}{2}}\leq
\eta_{j-1}^{\frac{3}{4}}$,  such that the transformation $\bar
\theta_{j-1}=\bar \theta_j+h_j(\bar \theta_j,\varphi)$ conjugates
the system
\begin{eqnarray} \label{sys-18}\left\{
\begin{array}{ll} & \dot{\bar{\theta}}_{j-1} = \rho_{f} +\overline{g}_{j-1}(\varphi) +
\overline{f}_{j-1}(\bar \theta_{j-1},\varphi) \\ & \dot{\varphi} =
\omega=(1,\alpha)
\end{array} \right.
\end{eqnarray}
to
\begin{eqnarray*} \left\{
\begin{array}{ll} & \dot{\bar\theta}_j = \rho_{f} +\overline{g}_{j}(\varphi) +
\overline{f}_{j}(\bar \theta_{j},\varphi) \\ & \dot{\varphi} =
\omega=(1,\alpha)
\end{array} \right.
\end{eqnarray*}
with $\overline{g}_j=\overline{g}_{j-1}+[\overline{f}_{j-1}(\bar
\theta_j,\varphi)]_{\bar \theta_j},$ where
$[\overline{f}_j]_{\bar\theta_j}$ means the mean value of
$\overline{f}_j$ in $\bar\theta_j$ over $\mathbb{T}$.

When $j=\nu$, under the transformation
$\bar{\theta}_{\nu-1}=\bar\theta_{\nu}+h_{\nu}(\bar
\theta_\nu,\varphi)$ (mod 1), the fibred equation of
$(\ref{sys-18})$ becomes
\begin{equation}\label{lemma5-3-1}
\dot{\bar\theta}_\nu=\rho_{f} +\overline{g}_{\nu-1}(\varphi)+
\overline{f}_{\nu-1}(\bar\theta_\nu,\varphi)-\partial_{\omega}h_{\nu}-(\rho_{f}
+\overline{g}_{\nu-1}(\varphi))\frac{\partial h_{\nu}}{\partial
\bar\theta_\nu}+h.o.t
\end{equation}where $h.o.t$ means higher order terms of $\overline{f}_{\nu-1}$ and $h_\nu$. So the homological equation is
\begin{equation}\label{coho-equ-1}
\partial_{\omega}h_{\nu}+(\rho_{f}
+\overline{g}_{\nu-1}(\varphi))\frac{\partial h_{\nu}}{\partial
\bar\theta_\nu}=\overline{f}_{\nu-1}(\bar\theta_\nu,\varphi)-[\overline{f}_{\nu-1}(\bar\theta_\nu,\varphi)]_{\bar\theta\nu},
\end{equation}

By the definition of $\bar g_{j}, j=1,2,\cdots,\nu$, we have $\|\bar
g_\nu\|_{\tilde r_{\nu-1}}\leq \|\bar g\|_{\tilde
r_0}+\sum_{j=0}^{\nu-1}\eta_j<2\varepsilon_{n-1}^{1/2}$. Moreover, by our
selection of $\varepsilon_0$ and the definition of
$\varepsilon_{n-1}, \widetilde{r}_{\nu-1},\widetilde{s}_{\nu-1},
\eta_{\nu-1}$, we have
\begin{eqnarray*}
K_\nu&=&\left[\frac{1}{\sigma_{\nu}}\ln\frac{1}{\eta_{\nu-1}}\right]\leq
\frac{4\cdot 3^{\nu-1}Q_n^3}{r_0}\ln\frac{1}{2\varepsilon_{n-1}}\\
&\leq & \frac{8Q_n^3\cdot 3^{2^nc_1\tau U\ln
Q_n}}{r_0}\ln\frac{1}{2\varepsilon_{n-1}}\\&=&\frac{8Q_n^{\frac{2^nc\tau
U}{32(\tau+3)}+3}}{r_0}\ln\frac{1}{2\varepsilon_{n-1}}\\
&<&
 \varepsilon_{n-1}^{-\frac{1}{16(\tau+3)}-\frac{1}{12(2\tau+3)}}\\
&<&\left(\frac{\gamma^2}{2\varepsilon_{n-1}^{1/2}}\right)^{\frac{1}{2\tau+3}}=K^{(n-1)}
\end{eqnarray*}
for $\nu\leq N$. By the assumption that $\rho_f\in
DC_\omega(\gamma,\tau)$, the conditions of Proposition
\ref{Diagonal} are satisfied. We can apply this proposition,
obtaining an approximate solution of the homological equation
(\ref{coho-equ-1}) with estimates
\begin{displaymath}
\|h_{\nu}\|_{\widetilde{s}_{\nu-1}-\frac{\delta_{\nu}}{2},\widetilde{r}_{\nu-1}-\frac{\sigma_{\nu}}{2}}
<\eta_{\nu-1}^{\frac{3}{4}}<\frac{\delta_{1}}{2^{\nu}}=\delta_{\nu}/2,
\end{displaymath}
and the error satisfying
$$\|\tilde P_{\nu}\|_{\widetilde{s}_{\nu-1}-\frac{\delta_{\nu}}{2},\widetilde{r}_{\nu-1}-\frac{\sigma_{\nu}}{2}}
\leq \frac{4\eta_{\nu-1}^2}{\gamma (\sigma_\nu/2)^{3+\tau}}\leq
\frac{2\eta_{\nu-1}^{\frac{3}{2}}}{5}.$$ Now (\ref{lemma5-3-1})
becomes
\begin{equation}
\dot{\bar\theta}_\nu=\rho_f+\bar g_\nu(\varphi)+\bar f_\nu(\bar
\theta_\nu,\varphi),
\end{equation}
where
$\overline{g}_\nu=\overline{g}_{\nu-1}+[\overline{f}_{\nu-1}(\bar
\theta_\nu,\varphi)]_{\bar \theta_\nu}$, and
\begin{eqnarray*}
&& \overline{f}_{\nu}(\bar\theta_\nu,\varphi)+
\overline{f}_{\nu}(\bar\theta_\nu,\varphi)\frac{\partial
h_{\nu}}{\partial \bar\theta_\nu}
\\&=&\left(\overline{f}_{\nu-1}(\bar\theta_\nu+h_\nu(\bar\theta_\nu,\varphi),\varphi)-\overline{f}_{\nu-1}(\bar\theta_\nu,\varphi)\right)
+\widetilde{P}_{\nu}
-[\overline{f}_{\nu-1}(\bar\theta_\nu,\varphi)]_{\bar\theta_\nu}\frac{\partial
h_{\nu}}{\partial \bar\theta_\nu}
\end{eqnarray*}
By the mean value theorem and the Cauchy estimate, we have
\begin{eqnarray*}
\|\overline{f}_{\nu-1}(\bar\theta_\nu+h_\nu(\bar\theta_\nu,\varphi),\varphi)-
\overline{f}_{\nu-1}(\bar\theta_\nu,\varphi)\|_{\widetilde{s}_{\nu},\widetilde{r}_{\nu}}
\leq\frac{\eta_{\nu-1}}{\delta_{\nu}/2}
\|h_{\nu}\|_{\widetilde{s}_{\nu-1}-\frac{\delta_{\nu}}{2},\widetilde{r}_{\nu-1}-\frac{\sigma_{\nu}}{2}}
\leq \frac{\eta_{\nu-1}^{\frac{3}{2}}}{5}.
\end{eqnarray*}
It follows that
\begin{eqnarray*}
\|\overline{f}_{\nu}(\bar\theta_\nu,\varphi)\|_{\widetilde{s}_{\nu},\widetilde{r}_{\nu}}&\leq&(1+\eta_{\nu-1}^{\frac{3}{4}})
\Big(\frac{\eta_{\nu-1}^{\frac{3}{2}}}{5}+\frac{2\eta_{\nu-1}^{\frac{3}{2}}}{5}+\frac{\eta_{\nu-1}^{\frac{3}{2}}}{5}\Big)
\leq\eta_{\nu-1}^{\frac{3}{2}}=\eta_{\nu}.
\end{eqnarray*}

To finish the proof, we estimate the size of
$\overline{f}_{N}(\bar \theta_N,\varphi)$. By the choice  of $U=U(\a)$ (cf. Corollary \ref{cor2-1}), we have
$Q_n^{U}\geq \ln Q_{n+1}$. As a consequence,
$$
(\frac{3}{2})^{N}-1\geq (\frac{3}{2})^{2^nc_1\tau U\ln Q_n}-1\geq
Q_n^{ 2^{n-1}c_1\tau U\ln (3/2)}.
$$
Therefore,
\begin{eqnarray*}
\|\overline{f}_{N}\|_{\widetilde{s}_{N},\widetilde{r}_{N}}&\leq&
\eta^{(\frac{3}{2})^{N}}
=\eta e^{-((\frac{3}{2})^{N}-1)\ln\frac{1}{\eta}}\\
&\leq&\eta e^{-Q_n^{2^{n-1}c_1\tau U \ln (3/2)}2^nc\tau U}\\
&<&\eta\frac{1}{Q_{n+1}^{2^{n+1}c\tau U}}=\varepsilon_{n},
\end{eqnarray*}
which is possible since by our selection $\varepsilon_n
\leq\varepsilon_0 < e^{-2c\tau U}$.

For simplicity, we just denote $\bar \theta_N$ by $\bar\theta_+$. We
denote $H_{\nu}(\bar \theta_\nu,\varphi)=(\bar
\theta_\nu+h_{\nu}(\bar \theta_\nu,\varphi)\ \textrm{mod}\
1,\varphi)$ for $1\leq \nu\leq N$. Since all these transformations
are close to the identity, then by Lemma \ref{fib-rota-2}, we have
$\rho(\omega, \rho_{f} +\overline{g}_{\nu}(\varphi) +
\overline{f}_{\nu}(\bar \theta_\nu,\varphi) )=\rho_f$.  Let
$$\overline{H}_{\nu}(\bar\theta_\nu,\varphi)=H_{1}\circ \cdots \circ H_{\nu-1}\circ H_{\nu} (\bar\theta_\nu,\varphi).$$
Then $\overline{H}_{\nu}(\bar\theta_\nu,\varphi)$ is analytic in
$D(\widetilde{s}_{\nu},\widetilde{r}_{\nu})$ and
$$\|\partial_{\bar\theta_\nu}(\Pi_{1}\circ
\overline{H}_{\nu}(\bar\theta_\nu,\varphi))\|_{
\widetilde{s}_{\nu},\widetilde{r}_{\nu}  } \leq
\prod_{j=1}^{\nu}(1+\eta_{j}^{\frac{3}{4}}).$$ If we rewrite
$\overline{H}_{N}(\bar\theta_+,\varphi)=(\bar\theta_++\widetilde{h}(\bar\theta_+,\varphi)\
\textrm{mod}\ 1,\varphi)$, then we have
\begin{eqnarray*}
&&\|\widetilde{h}(\bar\theta_+,\varphi)\|_{\widetilde{s}_{N-1}-\frac{\delta_N}{2},\widetilde{r}_{N-1}-\frac{\sigma_N}{2}}\\
&=&\|\Pi_{1}\circ\overline{H}_{N}(\bar\theta_+,\varphi)-\bar\theta_+\|_
{\widetilde{s}_{N-1}-\frac{\delta_N}{2},\widetilde{r}_{N-1}-\frac{\sigma_N}{2}}\\
&\leq&
\|\Pi_{1}\circ(\overline{H}_{N}(\bar\theta_+,\varphi)-\overline{H}_{N-1}(\bar\theta_+,\varphi))\|
_{\widetilde{s}_{N-1}-\frac{\delta_N}{2},\widetilde{r}_{N-1}-\frac{\sigma_N}{2}}\\
&&{}+\cdots+\|\Pi_{1}\circ(\overline{H}_{2}(\bar\theta_+,\varphi)-\overline{H}_{1}(\bar\theta_+,\varphi))\|_
{\widetilde{s}_{1}-\frac{\delta_2}{2},\widetilde{r}_{1}-\frac{\sigma_2}{2}}\\
&&{}+\|\Pi_{1}\circ\overline{H}_{1}(\bar\theta_+,\varphi)-\bar\theta_+\|_
{\widetilde{s}_{0}-\frac{\delta_1}{2},\widetilde{r}_{0}-\frac{\sigma_1}{2}}\\
&\leq& \sum_{\nu=1}^{N}
\prod_{j=1}^{\nu}(1+\eta_{j}^{\frac{3}{4}})\|h_{\nu}\|_{s_{\nu},r_{\nu}}\leq
2\eta^{\frac{3}{4}}<4\varepsilon_{n-1}^{3/4}.
\end{eqnarray*}
Similarly, we have
\begin{eqnarray}
\|\frac{\partial\widetilde{h}}{\partial
\bar\theta_+}(\bar\theta_+,\varphi)\|_{\widetilde{s}_{N},\widetilde{r}_{N}}
\leq 4\varepsilon_{n-1}^{3/4},
\end{eqnarray}
where $\Pi_{1}:\R\times\R^2 \rightarrow \R$ denotes the natural
projection to the first variable. 

In conclusion, let
$\overline{g}_{+}(\varphi)=\overline{g}_{N}(\varphi),$
$\overline{f}_{+}(\bar\theta_{+},\varphi)=\overline{f}_{N}(\bar\theta_{+},\varphi)$.
Then
$\overline{H}(\bar\theta_+,\varphi)=(\bar\theta_{+}+\widetilde{h}(\bar\theta_{+},\varphi)\
\textrm{mod} \ 1,\varphi)$ transforms the system $(\ref{sys-2})$ to
\begin{eqnarray} \left\{
\begin{array}{ll} & \dot{\bar\theta}_{+} = \rho_{f} +\overline{g}_{+}(\varphi) +
\overline{f}_{+}(\bar\theta_{+},\varphi) \\ & \dot{\varphi} =
\omega=(1,\alpha)
\end{array} \right.
\end{eqnarray}
with estimate
\begin{equation} \label{differenceg}
\|\overline{g}_{+}(\varphi)-\overline{g}(\varphi)\|_{r_{n+}}\leq\sum_\nu\|\overline{f}_{\nu}\|
_{\widetilde{s}_\nu,\widetilde{r}_\nu}<\sum_\nu\eta_\nu<4\varepsilon_{n-1}.
\end{equation}
It follows that
$\rho_{f} +\overline{g}_{+}(\varphi) +
\overline{f}_{+}(\bar\theta_{+},\varphi)\in \tF_{s_{n+},r_{n+} }(\rho_f,
2\varepsilon_{n-1}^{1/2},\varepsilon_n)$.
\qed

\subsubsection{End of one step.}
In the first step, we eliminate the non-resonant terms of
$g(\varphi)$ and as a result the transformation we obtain  is not close to the
identity. In order to get rotations reducibility results, we  need
to inverse the first step, which means conjugating back by the transformation of the first step.

\begin{Lemma}\label{lemma4-3}
Under the assumptions of Lemma \ref{lemma4-2}, there exists
 $\widetilde{H}\in
C_{s_{n},r_{n}}^\omega(\T\times\T^2,\T\times\T^2)$ with estimates
$$\|\widetilde{H}-id\|_{s_n,r_n}\leq 4\varepsilon_{n-1}^{\frac{3}{4}},$$
$$\|D(\widetilde{H}-id)\|_{s_n,r_n}\leq 4\varepsilon_{n-1}^{\frac{3}{4}},$$
 such that
$\widetilde{H}$ conjugates the system $(\ref{sys-13})$ to
\begin{eqnarray}\label{sys-17} \left\{
\begin{array}{ll} & \dot{ {\theta}}_+ = \rho_f + {g}_+(\varphi) +
 {f}_+({ {\theta}}_+,\varphi) \\ & \dot{\varphi} =
\omega=(1,\alpha)
\end{array} \right.
\end{eqnarray}
with $  \rho_f + {g}_+(\varphi) +
 {f}_+({ {\theta}}_+,\varphi) \in \tF_{ s_{n},r_{n} }
(\rho_f, 4\widetilde\varepsilon_{n}, \varepsilon_{n} ) $.
\end{Lemma}
\Proof By Lemma \ref{lemma4-1}, there exists $h(\varphi)\in
C_{\overline{r}_n}^{\omega}(\T^2,\R)$ such that the system
$(\ref{sys-13})$ can be conjugated to $(\ref{sys-2})$. Set
$H(\overline{\theta},\varphi)=(\overline{\theta}+h(\varphi)\
\textrm{mod}\ 1,\varphi)$. Then we apply Lemma \ref{lemma4-2} to get
$\widetilde{h}\in C_{s_{n+},r_{n+}}^\omega(\T\times\T^2,\R)$ such
that
$\overline{H}(\bar\theta_+,\varphi)=(\bar\theta_{+}+\widetilde{h}(\bar\theta_{+},\varphi)\
\textrm{mod}\ 1,\varphi)$ conjugates the system $(\ref{sys-2})$ to
$(\ref{sys-12})$ without changing the fibred rotation number.

Let
$ {g}_+(\varphi)=\mathcal{T}_{Q_n}g(\varphi)+\overline{g}_{+}(\varphi)$,
$ {f}_+( {\theta}_+,\varphi)=\overline{f}_{+}( {\theta}_+-h(\varphi),\varphi)$,
then  under the transformation
$\theta_+=\overline{\theta}_+-h(\varphi)\ (\textrm{mod}\ 1)$, the
fibred equation of system $(\ref{sys-12})$ becomes
\begin{eqnarray*}
\dot{ {\theta}}_+&=&\rho_{f} +\overline{g}_{+}(\varphi) +
\overline{f}_{+}(\overline{\theta}_+-h(\varphi),\varphi)+\partial_{\omega}h(\varphi)\\
&=&\rho_f+ {g}_+(\varphi) +
 {f}_+( {\theta}_+,\varphi).
\end{eqnarray*}
Then we have
$$\| {g}_+(\varphi)-g(\varphi)\|_{r_n}\leq \|\overline{g}_{+}(\varphi)-\overline{g}(\varphi)\|_{r_{n+}}<
4\varepsilon_{n-1}.$$
Therefore,  we have
$\|{g}_+(\varphi)\|_{r_n}\leq
4\varepsilon_{n-1}+4\widetilde\varepsilon_{n-1}=4\widetilde\varepsilon_{n}.$
Using the same argument as in Lemma
\ref{lemma4-1}, we have $$ \rho_f+ {g}_+(\varphi) +
 {f}_+( {\theta}_+,\varphi)\in
\tF_{s_{n},r_{n}}(\rho_f, 4\widetilde\varepsilon_{n},
\varepsilon_{n}).$$

Let $\widetilde{H}=H\circ \overline{H}\circ H^{-1}$. Then
$$\Pi_1\circ\widetilde{H}( {\theta}_+,\varphi)= {\theta}_++\widetilde{h}( {\theta}_+
-h(\varphi),\varphi)\ \textrm{mod}\ 1,$$ and
$\widetilde{H}( {\theta}_+,\varphi)$ conjugates the system
$(\ref{sys-13})$ to $( \ref{sys-17})$.
Furthermore, the estimates
$$\|\widetilde{H}-id\|_{s_n,r_n}\leq
\|\widetilde{h}(\bar\theta_+,\varphi)\|_{s_{n+},r_{n+}}<
 4\varepsilon_{n-1}^{\frac{3}{4}},$$
$$\|D( \widetilde{H}-id)\|_{s_n,r_n}\leq \|D\widetilde{h}\|_
{s_{n+},r_{n+}}< 4\varepsilon_{n-1}^{\frac{3}{4}}$$
follows from (\ref{tran-1}) and (\ref{tran}). \qed

\subsection{Iteration and convergence.}

Summarizing conclusions of section \ref{sec-3}, we have the
following iteration lemma.

\begin{Lemma}\label{lem-main}
For any $\varepsilon_{0}>0,r_0>0,s_{0}>0,\gamma>0,\tau>2$,
$\omega=(1,\alpha)$ with $U=U(\alpha)<\infty$, $\varepsilon_n,
\widetilde\varepsilon_n,r_n,s_n$ are defined as in (\ref{seque}).
 Suppose   $\varepsilon_{0}$ is small enough such that
it satisfies (\ref{enseq}).
Then the following holds for all $n\geq 1$: If the system
\begin{eqnarray}\label{sys-5} \left\{
\begin{array}{ll} & \dot{\theta} = \rho_f +g_{n}(\varphi) +
f_{n}(\theta,\varphi) \\ & \dot{\varphi} = \omega=(1,\alpha)
\end{array} \right.
\end{eqnarray}
satisfies $\rho_f +g_{n}(\varphi) + f_{n}(\theta,\varphi)\in \tF_{
s_{n},r_{n} } (\rho_f, 4\widetilde\varepsilon_{n},\varepsilon_{n}
)$, then
 there exists $H_{n}: \T\times\T^2\rightarrow \T\times\T^2$  with estimates
\begin{itemize}
\item $\|H_{n}-id\|_{s_{n+1},r_{n+1}}\leq 4\varepsilon_{n}^{\frac{3}{4}}$
\item  $\|D( H_{n} -id)\|_{s_{n+1},r_{n+1}}\leq
4\varepsilon_{n}^{\frac{3}{4}}$.
 \end{itemize}
such that it  transforms the system $(\ref{sys-5})$ to
\begin{equation}\label{sys-14} \left\{
\begin{array}{ll} & \dot{\theta} = \rho_f +g_{n+1}(\varphi) +
f_{n+1}(\theta,\varphi) \\ & \dot{\varphi} = \omega=(1,\alpha)
\end{array} \right.
\end{equation}
with $\rho_f +g_{n+1}(\varphi) + f_{n+1}(\theta,\varphi)\in \tF_{
s_{n+1},r_{n+1} } (\rho_f,
4\widetilde\varepsilon_{n+1},\varepsilon_{n+1} )$.
\end{Lemma}

With carefully checking of the proof, more precisely if we only do the first
two steps as in section \ref{sec-4-2-1} and \ref{sec-4-2-2},
 then   the following iteration lemma holds; it will be
the basis for the proof of almost reducibility.

\begin{Lemma}\label{lem-main-al}
Under the assumptions of Lemma \ref{lem-main},  the following holds
for all $n\geq 1$: If the system
\begin{eqnarray}\label{sys-5+} \left\{
\begin{array}{ll} & \dot{\theta} = \rho_f +g_{n}(\varphi) +
f_{n}(\theta,\varphi) \\ & \dot{\varphi} = \omega=(1,\alpha)
\end{array} \right.
\end{eqnarray}
satisfy $\rho_f +g_{n}(\varphi) + f_{n}(\theta,\varphi)\in \tF_{
s_{n},r_{n} } (\rho_f, 2\varepsilon_{n-1}^{1/2},\varepsilon_{n} )$, then
 there exists $\tilde{H}_{n}: \T\times\T^2\rightarrow \T\times\T^2$
with estimate  $\|\tilde{H}_{n}-id\|_{s_{n+1},r_{n+1}}\leq
Q_{n+1}^2\varepsilon_{0} $ such that it  transforms the system
$(\ref{sys-5+})$ to
\begin{equation}\label{sys-16} \left\{
\begin{array}{ll} & \dot{\theta} = \rho_f +g_{n+1}(\varphi) +
f_{n+1}(\theta,\varphi) \\ & \dot{\varphi} = \omega=(1,\alpha)
\end{array} \right.
\end{equation}
with $\rho_f +g_{n+1}(\varphi) + f_{n+1}(\theta,\varphi)\in \tF_{
s_{n+1},r_{n+1} } (\rho_f, 2\varepsilon_{n}^{1/2},\varepsilon_{n+1} )$.
\end{Lemma}

\section{Proof of Theorem \ref{mainthm}}

In this section, we first use Lemma \ref{lem-main} and Lemma
\ref{lem-main-al} to prove Theorem \ref{mainthm} (a). Then as
corollaries, we prove the local denseness of linearization (Theorem
\ref{mainthm} (b))  and local denseness of mode-locking (Corollary
\ref{dense-mode}).

\subsection{Proof of Theorem \ref{mainthm}(a):}
Select  $\varepsilon_0$ such that it satisfies (\ref{enseq}) with $r_0=r, s_0=s$.
We first prove that if $\rho(\omega, \tilde{\rho}+
f(\theta,\varphi)) =\rho_f\in DC_\omega(\gamma,\tau)$ and $\|f\|_{s,r}\leq
\varepsilon_0/2$, then the system $(\omega,
\tilde{\rho}+f(\theta,\varphi))$ is $C^\infty$ rotations
linearizable.

Without loss of generality, we can rewrite  the system $(\omega,
\tilde{\rho}+ f(\theta,\varphi))$ as
\begin{eqnarray} \label{sys-15}\left\{
\begin{array}{ll} & \dot{\theta} = \rho_f +\widetilde{g}(\varphi) +
\widetilde{f}(\theta,\varphi) \\ & \dot{\varphi} = \omega=(1,\alpha)
\end{array} \right.
\end{eqnarray}
with  $\rho_f +\widetilde{g}(\varphi) +
\widetilde{f}(\theta,\varphi)\in
\tF_{s,r}(\rho_f,\varepsilon_{0},\varepsilon_{0} ).$ Since
$\rho_f\in DC_\omega(\gamma,\tau)$ and $\varepsilon_0$ satisfies the
inequality $(\ref{enseq})$, we can apply Lemma~\ref{lemma4-2} (in
fact without using Lemma \ref{lemma4-1}  and Lemma \ref{lemma4-3}),
and thus get $H_0 \in C^\omega_{s_1,r_1}(\T\times\T^2,
\T\times\T^2)$ which conjugates $(\ref{sys-15})$ to
\begin{eqnarray} \left\{
\begin{array}{ll} & \dot{\theta} = \rho_f +g_{1}(\varphi) +
f_{1}(\theta,\varphi) \\ & \dot{\varphi} = \omega=(1,\alpha)
\end{array} \right.
\end{eqnarray}
 with $\rho_f +g_{1}(\varphi) +
f_{1}(\theta,\varphi)\in
\tF_{s_{1},r_{1}}(\rho_f,4\widetilde\varepsilon_{1},\varepsilon_{1}
)$. Then we apply Lemma \ref{lem-main}, and inductively  we obtain sequence $H_i\in C^\omega_{s_{i+1},r_{i+1}}(\T\times\T^2,
\T\times\T^2)$, $i=1,\cdots,n-1$ such that $H^{(n)}=H_0\circ H_1\circ \cdots \circ
H_{n-1}$ conjugates the system $(\ref{qpf-foced})$ to $(\omega,
\rho_f +g_{n}(\varphi) + f_{n}(\theta,\varphi))$ with $ \rho_f
+g_{n}(\varphi) + f_{n}(\theta,\varphi) \in \tF_{s_n,r_n}(\rho_f,
4\widetilde\varepsilon_{n},\varepsilon_{n} )$.

Let $H=\lim_{n\rightarrow\infty}H^{(n)},$
$g_\infty=\lim_{n\rightarrow \infty}g_n$. Then the system
$(\omega,\tilde{\rho}+ f(\theta,\varphi))$ is conjugated by $H$ to
$(\omega,\rho_f+g_\infty(\varphi))$. The remaining task is to prove
that the transformation $H$ is  actually in $C^\infty$. Since
$$\|DH^{(n)}\|_{s_n,r_n}\leq \|DH_0\|_{s_1,r_1}\|DH_1\|_{s_2,r_2}\cdots \|DH_{n-1}\|_{s_n,r_n}\leq \prod_{i=0}^{n-1}(1+4\varepsilon_i^{3/4})<2,$$
then we get
$$\|H^{(n+1)}-H^{(n)}\|_{s_{n+1},r_{n+1}}\leq \|DH^{(n)}\|_{s_n,r_n}\|H_{n}-id\|_{s_{n+1},r_{n+1}}
\leq 8\varepsilon_n^{3/4}.$$

 By the definition of
$(\varepsilon_n)_{n\in\N}$, we know that for any $j\in\Z_+^3$, there
exists some $N\in\N$, so that for any $n\geq N$, we have
$Q_n^{4|j|}<\varepsilon_{n-1}^{-1/4}$, that is
$$Q_n^{4|j|}\varepsilon_{n-1}^{3/4}<\varepsilon_{n-1}^{1/2},\ \ \forall n\geq N.$$
Then by Cauchy estimates, if we denote $x:=(\theta,\varphi)\in
\T^3$, we have
$$|\frac{\partial^{|j|}}{\partial x^j}(H^{(n+1)}-H^{(n)})|\leq r_{n+1}^{-|j|}\|H^{(n+1)}-H^{(n)}\|_{s_{n+1},r_{n+1}}
\leq Q_{n+1}^{4|j|}\varepsilon_n^{3/4}<\varepsilon_n^{1/2} $$ for
any $n\geq N-1$. This guarantees the limit
$H=\lim_{n\rightarrow\infty}H^{(n)}$ belongs to $C^\infty$. As a
consequence, we have $g_\infty\in C^\infty(\T^2,\R)$.

The proof of the $C^\infty$ almost reducibility is exactly the same as the proof
of the $C^\infty$ rotations reducibility.  The only difference is that
we use  Lemma \ref{lem-main-al} instead of Lemma \ref{lem-main}.
Note in this case, if
$\widetilde{H}=\lim_{n\rightarrow\infty}\widetilde
H_0\circ\cdots\circ \widetilde H_{n-1},$ then $\widetilde{H}$
actually diverges.\qed

\subsection{Proof of Theorem \ref{mainthm} (b):}

Suppose that  $\varepsilon_{0}$ is selected according to $(\ref{enseq})$ with $r_0=r, s_0=s$.
 Consider the system
$(\omega, \tilde{\rho}+ f(\theta,\varphi))$ with $ \rho(\omega,
\tilde{\rho}+ f(\theta,\varphi)) =  \rho_f\in DC_\omega(\gamma,\tau)$
and  $\|f\|_{s,r}<\varepsilon_0/2$. For any $\epsilon>0$, there
exists $N\in \N$ such that $\frac{4}{Q_{N-1}^{2^{N-1}c\tau
U}}<\epsilon$. Then by Lemma $\ref{lem-main-al}$, there exists a map
$H\in C^\omega(\T\times\T^2,\T\times\T^2)$ such that the system
$(\omega, \tilde{\rho}+ f(\theta,\varphi))$ can be conjugated to
$(\omega,\rho_f+\overline{g}(\varphi)+\overline{f}(\theta,\varphi))$
with
$$\|\overline{g}\|_{\overline{r}}<\epsilon^{1/2},\quad
\|\overline{f}\|_{\overline{s},\overline{r}}<\epsilon\frac{1}{Q_{N}^{2^Nc\tau
U}}.$$
$$\|\Pi_1\circ
H-id\|_{\overline{s},\overline{r}}<Q_{N-1}^2 \varepsilon_{0}, \quad
\|\frac{\partial}{\partial\theta}(\Pi_1\circ
H-id)\|_{\overline{s},\overline{r}}<Q_{N-1}^2 \varepsilon_{0}.$$ where
$\overline{s}=s_{N-1}$, $\overline{r}=r_{N-1}$. Rewrite
$(\omega,\rho_f+\overline{g}(\varphi)+\overline{f}(\theta,\varphi))$
as
$(\omega,\widetilde{\overline{\rho}}+\widetilde{\overline{g}}(\varphi)+\widetilde{\overline{f}}(\theta,\varphi))$,
where
$\widetilde{\overline{\rho}}=\rho_f+\widehat{\overline{g}}(0)$,
$\widetilde{\overline{g}}(\varphi)=\mathcal{T}_{Q_N}\overline{g}(\varphi)$,
and
$\widetilde{\overline{f}}(\theta,\varphi)=\overline{f}(\theta,\varphi)+\mathcal{R}_{Q_N}g(\varphi)$.
Then let $\widetilde{s}=\overline{s}$,
$\widetilde{r}=\frac{\overline{r}}{2}$, we have
\begin{equation}\label{pert}
\|\widetilde{\overline{f}}\|_{\widetilde{s},\widetilde{r}}
\leq\|\overline{f}\|_{\overline{s},\overline{r}}+
\|\mathcal{R}_{Q_N}g(\varphi)\|_{\overline{s},\frac{\overline{r}}{2}}
<2\epsilon\frac{1}{Q_{N}^{2^Nc\tau U}}.
\end{equation}
We consider the reference system
$(\omega,\widetilde{\overline{\rho}}+\widetilde{\overline{g}}(\varphi))$:
it must be $C^{\omega}$ linearizable. The reason is that since
$\widetilde{\overline{g}}(\varphi)$ is a trigonometric polynomial,
if we let
$\partial_{\omega}\widetilde{h}(\varphi)=\widetilde{\overline{g}}(\varphi)$,
then  $\widetilde{h}(\varphi)$ is analytic.

If we conjugate back both
$(\omega,\widetilde{\overline{\rho}}+\widetilde{\overline{g}}(\varphi))$
and
$(\omega,\widetilde{\overline{\rho}}+\widetilde{\overline{g}}(\varphi)+\widetilde{\overline{f}}(\theta,\varphi))$
by the above transformation $H$, we obtain
$(\omega,\widetilde{f}(\theta,\varphi))$ and
$(\omega,\tilde{\rho}+f(\theta,\varphi))$. Here,
$(\omega,\widetilde{f}(\theta,\varphi))$ is clearly $C^\omega$
linearizable, since linearization is conjugacy invariant. As a
consequence of $(\ref{pert})$, we have
\begin{eqnarray*}
\|\tilde{\rho}+ f(\theta,\varphi)- \widetilde{f}(\theta,\varphi)\|_{\widetilde{s},\widetilde{r}}
&\leq& \|\frac{\partial}{\partial\theta}\Pi_1\circ
H\|_{\widetilde{s},\widetilde{r}}\cdot
\|\widetilde{\overline{f}}\|_{\widetilde{s},\widetilde{r}}<\epsilon.
\end{eqnarray*}
That is to say, $(\omega,\tilde{\rho}+ f(\theta,\varphi))$ is
$C^\infty$-accumulated
by the qpf circle flow
$(\omega,\widetilde{f}(\theta,\varphi))$, while
$(\omega,\widetilde{f}(\theta,\varphi))$ is $C^{\omega}$
linearizable.\qed

\subsection{Proof of Corollary \ref{dense-mode}}

To complete the picture in the local situation, we give the
following proposition:
\begin{Proposition}\label{rig-con}
Any rigid rotation flow $(\omega,\rho)$ can be $C^\infty$ accumulated by flows $(\omega,\widehat{f}(\theta,\varphi))$, where
$(\omega, \widehat{f}(\theta,\varphi))$ is mode-locked.
\end{Proposition}
\Proof For any $\rho >0,$ $\varepsilon>0$, there exists $k\in\Z^2$
such that $|\la k,\omega\ra-\rho| <\frac{\varepsilon}{2}$, since
$\omega$ is rationally independent. Under the transformation
$\theta=\overline{\theta}+\la k,\varphi\ra\ (\textrm{mod}\ 1)$,
$(\omega,\la k,\omega\ra)$ is conjugated to $(\omega,0)$. Now we
choose a qpf circle flow
$(\omega,\frac{\varepsilon}{4\pi}\sin4\pi\overline{\theta})$, which
can be viewed as projective action of a uniformly hyperbolic
$sl(2,\R)$ quasi-periodic flow:
\begin{eqnarray*}
\left\{ \begin{array}{l}\dot{x}=\left( \begin{array}{ccc}
 \frac{\varepsilon}{2} &  0\cr
 0 &   -\frac{\varepsilon}{2}\end{array} \right)x \\
\dot{\varphi}=\omega
\end{array}. \right.
\end{eqnarray*}
By Proposition  \ref{mode}, $(\omega,\frac{\varepsilon}{4\pi}
\sin4\pi\overline{\theta})$ is mode-locked. Invert the
transformation $\theta=\overline{\theta}+\la k,\varphi\ra$ (mod 1)
both for $(\omega,0)$ and $(\omega,\frac{\varepsilon}{4\pi}
\sin4\pi\overline{\theta})$, we get $(\omega,\la k,\omega\ra)$ and
$(\omega,\widehat{f})$. A direct calculation shows that
$\|\widehat{f}(\theta,\varphi)-\rho\|_{s,r}<\frac{\varepsilon}{2},$  while
$(\omega,\widehat{f})$ is mode-locked.\qed

Now we can finish the proof of Corollary \ref{dense-mode}. For any
qpf circle flow $(\omega, \tilde{\rho}+f(\theta,\varphi))$ with
$\rho(\omega, \tilde{\rho}+f(\theta,\varphi))=\rho_f\in
DC_\omega(\gamma,\tau)$ and $\|f\|_{r,s}<\varepsilon_0/2$, first we
perturb $(\omega, \tilde{\rho}+f(\theta,\varphi))$ to $C^\omega$
linearizable qpf circle flow
$(\omega,\widetilde{f}(\theta,\varphi))$ by Theorem \ref{mainthm}
(b). Then $(\omega,\widetilde{f}(\theta,\varphi))$ is $C^{\infty}$-accumulated 
by
mode-locked qpf circle flow by Proposition \ref{rig-con}.\qed

\section*{Acknowledgements}

R.Krikorian was supported by the ``Chaire d'Excellence''  Labex MME-DII at the University of Cergy-Pontoise and by the ANR project BEKAM (ANR-15-CE40-0001).    J. Wang has been supported by a research fellowship of the Alexander-
Humboldt-Foundation, NNSF of China (11601230) and  Natural Science Foundation of Jiangsu Province, China(BK20160816).  J. You was supported by NNSF of China (11471155) and 973 projects of China (2014CB340701).
Q. Zhou was  supported by  Fondation Sciences Math\'{e}matiques de Paris (FSMP) and NNSF of China (11671192).


\begin{thebibliography}{99}




\bibitem {AR}
V.I. Arnold,  
Small denominators I. On the mapping of a circle into itself.
Akad.Nauk.Math,\textbf{25} (1961),21-86.



\bibitem{Aglobal}
 A. Avila,
Global theory of one-frequency Schr\"{o}dinger operators, Acta Math. \textbf{215}, 1-54 (2015).


\bibitem {Aac}
 A. Avila,
Almost reducibility and absolute continuity,
http://w3.impa.br/~avila/

\bibitem{A2}
 A. Avila, KAM, Lyapunov exponents and the spectral dichotomy for one-frequency schr\"odinger operators. In preparation.



\bibitem{AFK}
\newblock A. Avila, B. Fayad and  R. Krikorian,
 \newblock  \ A KAM scheme for $\SL(2,\R)$ cocycles with Liouvillean
 frequencies.
\newblock Geom. Funct. Anal. \textbf{21} (2011), 1001-1019.

\bibitem{AJ05}
\newblock A. Avila  and  S. Jitomirskaya, 
 \newblock The ten Martini problem,
\newblock Annals of Mathematics,\textbf{170} (2009)  303-342.

\bibitem{AK06}
\newblock A. Avila and R. Krikorian
\newblock Reducibility or non-uniform hyperbolicity
for quasiperiodic Schr\"{o}dinger cocycles.
\newblock Annals of Mathematics. \textbf{164} (2006), 911-940.


\bibitem{bellissard/simon:1982}
J.~B\'{e}llissard and B.~Simon.
\newblock Cantor spectrum for the almost {Mathieu} equation.
\newblock J. Funct. Anal., \textbf{48 (3)}  (1982), 408--419.

\bibitem{BJ}
\newblock K. Bjerkl\"{o}v and T. J\"{a}ger, 
\newblock Rotation numbers for quasiperiodically forced circle maps -
mode-locking vs strict monotonicity.
\newblock  J. Am. Math.
Soc. \textbf{22} (2009), 353-362.



\bibitem{DS75}
\newblock E. Dinaburg and Ya. Sinai, 
\newblock The one-dimensional
Schr\"{o}dinger equation with a quasi-periodic potential.
\newblock Funct. Anal. Appl. \textbf{9} (1975), 279-289.


\bibitem{DGO89}
\newblock M. Ding, C. Grebogi, E. Ott, 
 \newblock Evolution of attractors in quasiperiodically forced systems: From quasiperiodic to strange nonchaotic to chaotic,
\newblock Physical Review A, 39(\textbf{5}) (1989)  2593-2598.


\bibitem{E92}
\newblock H. Eliasson, 
\newblock Floquet solutions for the one-dimensional quasiperiodic Schr\"{o}dinger
equation.
\newblock Comm.Math.Phys. \textbf{146} (1992), 447-482.

\bibitem{FK}
\newblock B. Fayad and R. Krikorian,
\newblock \ Rigidity results for quasiperiodic
$\SL(2,\R)$-cocycles.
\newblock  Journal of Modern Dynamics. \textbf{3}, (2009), no. 4.
479-510.





\bibitem{He79}
\newblock M. Herman, 
\newblock Sur la conjugaison diff\'rentiable des diff\'eomorphismes du cercle \`a des rotations,
\newblock  Inst. Hautes Etudes Sci. Publ. Math. (1979) 5-233.

\bibitem{He83}

\newblock M. Herman, 
\newblock Une {m\'{e}thode} pour minorer les exposants de {{L}yapunov} et
  quelques exemples montrant le {caract\`{e}re} local d'un {th\'{e}or\`{e}me}
  {d'Arnold} et de {Moser} sur le tore de dimension 2,
\newblock  Comment.\ Math.\ Helv. \textbf{58} (1983) 453-502.

\bibitem{HoY}
\newblock X. Hou and J. You, 
\newblock Almost reducibility and non-perturbative reducibility of quasiperiodic linear
systems.
\newblock  Invent. Math.  \textbf{190} (2012) 209-260.






\bibitem{JS06}
\newblock T. J\"{a}ger, J. Stark, 
\newblock Towards a classification for quasiperiodically forced circle homeomorphisms,
\newblock  Journal of the LMS, \textbf{73(3)} (2006) 727-744.



\bibitem{Kr}
\newblock R. Krikorian,
\newblock Reducibility, differentiable
rigidity and Lyapunov exponents for quasiperiodic cocycles on $\T \times SL(2,\R)$,
\newblock http://arxiv.org/abs/math/0402333.

\bibitem{K15}
\newblock R. Krikorian,
\newblock On almost reducibility of circle diffeomorphisms,
\newblock in preparation.




\bibitem{puig}
\newblock J. Puig, 
\newblock Cantor spectrum for the Almost Mathieu operator.
\newblock Comm.Math.Phys. \textbf{244} (2004), 297-309.


\bibitem{Ji12}
\newblock J. Wang, 
\newblock Lower dimensional invariant tori for quasiperiodically forced circle diffeomorphisms,
\newblock  J. Diff. Equations  \textbf{253} (2012) 1489-1543.


\bibitem{Yo84}
\newblock J. C. Yoccoz, 
\newblock Conjugaison diff\'erentiable des diff\'eomorphismes du cercle dont le nombre de rotation verifie une condition diophantienne,
\newblock  Ann. Sci. Ecole Norm. Sup. \textbf{17} (1984) 333-359.


\bibitem{Yo02}
\newblock J. C. Yoccoz, 
\newblock Analytic linearization of circle diffeomorphisms,
\newblock  in Dynamical Systems and Small divisors(Cetraro, 1998),
Lecture Notes in Math. \textbf{1784}, Springer-Verlag, New York,
(2002)  125-173.




\bibitem{YZ1}
 J. You and Q. Zhou, 
 \newblock Embedding
 of analytic quasi-periodic cocycles into analytic quasi-periodic linear systems and its
 applications.
 \newblock Comm. Math. Phys.  \textbf{323}(2013), 975-1005.




\bibitem{ZW}
\newblock Q.Zhou and J.Wang, 
\newblock  Reducibility results for quasiperiodic cocycles  with Liouvillean
frequency. J. Dyn. Diff. Equat., \textbf{24} (2012) 61-83.
\newblock



\end{thebibliography}
\end{document}